\documentclass[11pt]{article}
\usepackage{amsmath,amssymb,amsthm,latexsym}
\usepackage[dvips]{graphics}
\usepackage{epsfig, rotate}
\usepackage[english]{babel}
\usepackage{amsfonts}
\usepackage{array}
\usepackage{color}

\newtheorem{theorem}{{Theorem}}

\numberwithin{equation}{section} \numberwithin{theorem}{section}
\numberwithin{lemma}{section}\numberwithin{corollary}{section}

\usepackage%[pctex32]
{graphics}
\usepackage%[dvips]
{graphicx}
\numberwithin{equation}{section}

\newcommand{\bd}{\bold}
\def\td{\stackrel{d}{\rightarrow}}

\begin{document}

\noindent {\bf \Large  On Tests for Complete Independence of Normal Random Vectors}

\vspace{20pt}
\noindent{\bf Shuhua Chang$^1$, Yongcheng Qi$^2$}

\vspace{10pt}

{\small
\noindent $^1$Coordinated Innovation Center for Computable Modeling in Management Science, Tianjin University of Finance and Economics, Tianjin 300222, PR China. \\ Email: szhang@tjufe.edu.cn

\vspace{10pt}

\noindent $^2$Department of Mathematics and Statistics, University of Minnesota Duluth,
1117 University Drive, Duluth, MN 55812, USA.
\\Email: yqi@d.umn.edu

\date{\today}

\vspace{20pt}

\noindent{\bf Abstract.}
Consider a random sample of $n$ independently and identically distributed $p$-dimensional normal random vectors. A test statistic for
complete independence of high-dimensional normal distributions, proposed by Schott (2005),
is defined as the sum of squared Pearson's correlation coefficients. A modified test statistic has been proposed by Mao (2014). Under the assumption of complete independence, both test statistics are asymptotically normal if the limit $\lim_{n\to\infty}p/n$ exists and is finite.   In this paper,  we investigate the limiting distributions for both Schott's and  Mao's test statistics. We show that both test statistics, after suitably normalized, converge in distribution to the standard normal as long as both $n$
and $p$ tend to infinity. Furthermore, we show that the distribution functions of the test statistics can be approximated very well
by a chi-square distribution function with $p(p-1)/2$ degrees of freedom as $n$ tends to infinity regardless of how $p$ changes with $n$.

\vspace{20pt}

}

%\vspace{20pt}

\newpage

\section{Introduction}

In classical multivariate analysis, statistical methods have been developed mainly for data from designed experiments and
dimensions of the data are fixed or very small compared with the sample size.
Nowadays, new technology has generated various types of high-dimensional data sets such as financial data, consumer data,
modern manufacturing data, multimedia data, hyperspectral image data, internet data, microarray and DNA data.
A common feature for all these datasets is that their dimensions can be very large compared with
their sample sizes. See, e.g., Schott (2001, 2005, 2007), Ledoit and Wolf (2002), Fan, Peng and Huang (2005),
Bai et al. (2009), Chen et al (2010), Chen and Qin (2010), Fujikoshi et al. (2010), B\"uhlmann and van de Geer (2011),
Jiang et al (2012), Srivastava and Reid (2012).

Throughout the paper, $N_p(\bd{\mu}, \bd{\Sigma})$ denotes the $p$-dimensional normal distribution with mean vector $\bd{\mu}$ and covariance matrix $\bd{\Sigma}$, and $\bd{I}_p$ denotes the $p\times p$ identity matrix.  We assume that $\bd\Sigma$ is positive definite.
Write $\bd\Sigma=(\sigma(i,j))_{1\le i,j\le p}$. Then, $\bd\Gamma=(\rho_{ij})_{1\le i,j\le p}$ is the correlation matrix of $\bd\Sigma$ given by $\rho_{ij}=\sigma(i,j)/\sqrt{\sigma(i,i)\sigma(j,j)}$.

Assume that a $p$-dimensional random vector $\mathbf{x}=(x_1, \cdots, x_p)'$ has a distribution $N_p(\mathbf{\mu}, \mathbf{\Sigma})$. We are interested in
testing whether  the $p$ components $x_1, x_2, \cdots, x_p$ are independent or equivalently testing whether  the covariance matrix
$\mathbf{\Sigma}$ is diagonal.  Then,
the test can be written as
\begin{equation}\label{approach}
H_0: \bd\Gamma=\mathbf{I}_p\ \ \mbox{vs}\ \ H_a: \bd\Gamma \ne \mathbf{I}_p.
\end{equation}
In literature, \eqref{approach} is known as the test of complete independence.

Let $\mathbf{x}_1, \cdots, \mathbf{x}_n$ be i.i.d. from $N_p(\mu, \mathbf{\Sigma})$. Write
\[
\mathbf{x}_k=(x_{k1}, \cdots, x_{kp})', ~~~~k=1, \cdots, n.
\]
Define
\begin{equation}\label{tube}
r_{ij}=\frac{\sum_{k=1}^n(x_{ki}-\bar{x}_i)(x_{kj} - \bar{x}_j)}
{\sqrt{\sum_{k=1}^n(x_{ki}-\bar{x}_i)^2\cdot \sum_{k=1}^n(x_{kj}-\bar{x}_j)^2}},
\end{equation}
where $\bar{x}_i=\frac{1}{n}\sum_{k=1}^nx_{ki}$ and $\bar{x}_j=\frac{1}{n}\sum_{k=1}^nx_{kj}$. Then,
$\bd{R}_{n} :=(r_{ij})_{p\times p}$ is the sample correlation matrix based on the  $p$-dimensional random vectors  $\bd{x}_1, \cdots, \bd{x}_n$.

In classic multivariate analysis when $p$ is a fixed integer,  the likelihood  method is a nice approach to test \eqref{approach}.
From Bartlett (1954) or Morrison (2005), the likelihood ratio test statistic is a function of the determinant of $\mathbf{R}_n$.
When $p=p_n$ depends on $n$ and $p_n\to\infty$, the limiting distribution of the  likelihood ratio test statistic has been obtained
in Jiang and Yang (2013), Jiang, Bai and Zheng (2013) and Jiang
and Qi (2015), and  the likelihood ratio method can still be used to test \eqref{approach}.
However, the likelihood ratio method fails
when $p\ge n$, since the sample correlation matrix $\mathbf{R}_n$ is singular and  the corresponding test statistic is degenerate.
A natural requirement for non-singularity of $\bd{R}_{n}$ is $p<n$.

Schott (2005) considers the following test statistic
\[
t_{np}=\sum_{1\le j<i\le p}r_{ij}^2.
\]
Assume that the null hypothesis of \eqref{approach} holds and
$\lim_{n\to\infty}p/n=\gamma\in (0, \infty)$.
Schott (2005) proves that $t_{np}-\frac{p(p-1)}{2(n-1)}$ converges in distribution
to a normal distribution with mean $0$ and variance $\gamma^2$, that is,
\begin{equation}\label{schott}
t_{np}^*:=\frac{t_{np}-\frac{p(p-1)}{2(n-1)}}{\tau_{np}}\td N(0,1),
\end{equation}
where $\tau_{np}^2=\frac{p(p-1)(n-2)}{(n-1)^2(n+1)}$.

%An approximate  level $\alpha$ test would reject the null hypothesis
%in \eqref{approach} if $t_{np}-\frac{p(p-1)}{2(n-1)}>z_{\alpha}\sqrt{\frac{p(p-1)(n-2)}{(n-1)^2(n+1)}}$, where $z_{\alpha}$ is the $\alpha$ %level critical value of the
%standard normal distribution.

It is worth noting that the same test statistic $t_{np}$ is also proposed by Srivastava (2005).
%It is noted in Srivastava (2005) that the test statistic after normalized
%slightly differently converges to normality at a rather slow rate.
Srivastava (2005, 2006) also considers
a test statistic which is based on the Fisher's z-transformation and originally proposed by Chen and Mudholkar (1990):
\[
Q_{np}=\frac{(n-3)\sum_{1\le j<i\le p}z_{ij}^2-\frac12p(p-1)}{\sqrt{p(p-1)}},
\]
where $z_{ij}=\frac12\log\frac{1+r_{ij}}{1-r_{ij}}$.  From Srivastava (2005), such a test
has not been designed for large $p$. Instead, Srivastava (2005) proposes a test statistic $T_3$ which is related
to the sample covariances only. See Srivastava (2005, 2006) for  details. Under certain conditions, Srivastava (2005) shows that $T_3$
converges in distribution to the standard normal under the null hypothesis in \eqref{approach}.   A simulation study in
Srivastava (2006) indicates that $Q_{np}$ statistic is inferior as the test does not give a consistent nominal level when $n$ and $p$
are close.

Very recently, Mao (2014) proposes a new test for complete independence.
The new test statistic is closely related to Schott's test and is
defined by
\[
T_{np}=\sum_{1\le j<i\le p}\frac{r_{ij}^2}{1-r_{ij}^2}.
\]
It has been proved in Mao (2014) that $T_{np}$ is asymptotically normal
under the null hypothesis of \eqref{approach} and assumption that
$\lim_{n\to\infty}p/n=\gamma\in (0, \infty)$.

In this paper, we will remove the condition imposed on $p$ and assume only that $p=p_n\to\infty$ as $n\to\infty.$  We will show that
both $T_{np}$ and $t_{np}$ are asymptotically normal. We also establish a unified chi-square approximation for the distribution of $T_{np}$ and $t_{np}$ regardless of how $p$ changes with $n$.

The rest of the paper is organized as follows. The main results of the paper are given in section~\ref{main} and their proofs are postponed
until section~\ref{proof}.  A simulation study to compare the performance of several different approaches is reported in section~\ref{simulation}.

\section{Main Results}\label{main}

Our main results include three theorems. We first obtain the limiting distribution of the test statistic $T_{np}$ in a larger range
for $p$, and then establish a unified chi-square approximation for all $p\ge 2$. The limit distribution of $t_{np}$ is given in the third
theorem.

The first theorem states that Mao's (2014) test statistic $T_{np}$
is asymptotically normal as long as $p=p_n\to\infty$ as $n\to\infty$.

\begin{theorem}\label{thm1} Assume $p=p_n\to\infty$ as $n\to\infty$. Then, under the null hypothesis of \eqref{approach}
\begin{equation}\label{s2}
T_{np}^*:=\frac{T_{np}-\frac{p(p-1)}{2(n-4)}}{\sigma_{np}}\overset{d}\to N(0,1)
\end{equation}
as $n\to\infty$, where
\[
\sigma_{np}^2=\frac{p(p-1)(n-3)}{(n-4)^2(n-6)}.
\]
\end{theorem}

We expect that the limiting distribution of $(n-4)T_{np}$ is chi-squared with $p(p-1)/2$ degrees
of freedom when $p$ is fixed, and this will be confirmed in the following theorem.
For applications there seems a gap in the limiting
distributions of the test statistic $T_{np}$ as one has to distinguish whether $p=p_n$ converges or diverges.
Instead, under linear transformation we define a slightly different statistic as follows
\begin{equation}\label{TC}
T_{np}^c=\sqrt{p(p-1)}T^*_{np}+\frac{1}{2}p(p-1)=\sqrt{\frac{n-6}{n-3}}(n-4)T_{np}+\frac{1}{2}p(p-1)\left(1-\sqrt{\frac{n-6}{n-3}}\right).
\end{equation}
The statistic $T_{np}^c$ can fill this gap. Our second theorem reveals that
the chi-square distribution can be used to approach the distribution of $T_{np}^c$ no matter how $p=p_n$ changes with $n$.

\begin{theorem}\label{unified}
Let $p=p_n$, $n\ge 1$ be a sequence of positive integers with $p_n\ge 2$ for all large $n$. Then under the null hypothesis of \eqref{approach}
\begin{equation}\label{chisq}
\sup_x|P(T_{np}^c\le x)-P(\chi^2_{p(p-1)/2}\le x)|\to 0~~~\mbox{ as }n\to\infty.
\end{equation}
\end{theorem}

\vspace{10pt}

Theorem~\ref{unified} implies that  $T_{np}^c$ converges in distribution to a chi-square distribution with $p(p-1)/2$ degrees of freedom
uniformly over $p\ge 2$ as $n\to\infty$, that is,  the superium of the left-hand side of \eqref{chisq} over $p\ge 2$ converges
to zero as $n\to\infty$.

For comparison purpose, we extend Schott's statistic $t_{np}$ in the same manner. We will show that the central limit theorem
\eqref{schott} holds for all large $p$ and a chi-square approximation can also be applied to $t_{np}$ for small $p$.  Now we define
\begin{equation}\label{tc}
t_{np}^c=\sqrt{p(p-1)}t^*_{np}+\frac{1}{2}p(p-1)=\sqrt{\frac{n+1}{n-2}}(n-1)t_{np}+\frac{1}{2}p(p-1)(1-\sqrt{\frac{n+1}{n-2}}).
\end{equation}

\begin{theorem}\label{schott-all} (i) If $p=p_n\to\infty$ as $n\to\infty$, then \eqref{schott} holds under the null hypothesis of \eqref{approach}.\\
(ii) Let $p=p_n$ be a sequence of positive integers with $p_n\ge 2$ for all large $n$. Then, under the null hypothesis of \eqref{approach}
\begin{equation}\label{t-chisq}
\sup_x|P(t_{np}^c\le x)-P(\chi^2_{p(p-1)/2}\le x)|\to 0~~~\mbox{ as }n\to\infty.
\end{equation}
\end{theorem}

\vspace{10pt}

Assume $\alpha\in (0,1)$.
Let $z_\alpha$ and $\chi^2_\alpha(p(p-1)/2)$ denote the $\alpha$ level critical values for the standard normal distribution and
the chi-squared distribution with $p(p-1)/2$ degrees of freedom, respectively.

Based on \eqref{s2}, an approximate  level $\alpha$ test for \eqref{approach} has a critical region or rejection region
\begin{equation}\label{RT}
\mathcal{R}_{T}^*(\alpha)=\left\{T_{np}\ge \frac{p(p-1)}{2(n-4)}+ z_{\alpha}\sqrt{\frac{p(p-1)(n-3)}{(n-4)^2(n-6)}}\right\}.
\end{equation}
Based on the chi-square approximation \eqref{chisq}, an approximate  level $\alpha$ test rejects \eqref{approach} in the region
\begin{equation}\label{RTC}
\mathcal{R}_{T}^c(\alpha)=\left\{T_{np}\ge \frac{p(p-1)}{2(n-4)}\left(\sqrt{\frac{n-3}{n-6}}-1\right)+
\chi^2_{\alpha}(p(p-1)/2)\sqrt{\frac{(n-3)}{(n-4)^2(n-6)}}\right\}.
\end{equation}
While based on the normal approximation \eqref{schott} to Schott's test statistic $t_{np}$,
an approximate  level $\alpha$ test for \eqref{approach} has a rejection region
\begin{equation}\label{Rt}
\mathcal{R}_{t}^*(\alpha)=\left\{t_{np}\ge \frac{p(p-1)}{2(n-1)}+ z_{\alpha}\sqrt{\frac{p(p-1)(n-2)}{(n-1)^2(n+1)}}\right\}.
\end{equation}
Similarly, we have an approximate level $\alpha$ rejection region for test \eqref{approach}
\begin{equation}\label{RtC}
\mathcal{R}_{t}^c(\alpha)=\left\{t_{np}\ge \frac{p(p-1)}{2(n-1)}(\sqrt{\frac{n-2}{n+1}}-1)+
\chi^2_{\alpha}(p(p-1)/2)\sqrt{\frac{(n-2)}{(n-1)^2(n+1)}}\right\}
\end{equation}
 based on the chi-square approximation \eqref{chisq-all}.

\section{Simulation Study}\label{simulation}

Mao (2014) has conducted a simulation study and compared the performance of three test statistics including Mao's $T_{np}$, Schott's $t_{np}$ and Srivastava's $T_3$.   It has been reported in Mao (2014) that Mao's test statistic is comparable to
the other two test statistics in terms of the accuracy of sizes of the tests and outperforms in some models under weak dependence.

In this section we will carry out a finite-sample simulation study to compare the performance of Schott's (2005) $t_{np}$ and Mao's (2014) $T_{np}$ based  on the normal approximation and the chi-square approximation. We will not simply repeat
Mao's (2014) choices. Our focus is on the two test statistics $t_{np}$ and $T_{np}$  which are related to the sample correlations. More specifically,
we consider four normalized test statistics: $t^*_{np}$, $T^*_{np}$, $t^c_{np}$ and $T^c_{np}$. Their limiting distributions are determined by \eqref{schott},
\eqref{s2}, \eqref{t-chisq} and \eqref{chisq}, respectively, and the corresponding rejection regions for the four tests at level $\alpha$
are given by \eqref{Rt}, \eqref{RT}, \eqref{RtC} and \eqref{RTC}.

Let $\bd{\Sigma}_p^{(\rho)}$ denote a $p\times p$ matrix whose diagonal entries are equal to $1$ and all off-diagonal entries are equal
to $\rho$, where $\rho\in (-1, 1)$.   $\bd{\Sigma}_p^{(\rho)}$ is the covariance matrix of a normal random vector with all $p$ components
being standard normal random variables and covariances (and correlation coefficients) equal to $\rho$.  A random sample of size $n$ is drawn
from multivariate normal distribution $N_p(0, \bd{\Sigma}_p^{(\rho)})$ with the different choices for $n=15$, $30$, $60$, $100$, $200$,
$p=3$, $10$, $20$, $50$, $100$, $200$, and $\rho=0$, $0.02$. For each combination of
the choices on $n$, $p$ and $\rho$, the simulation experiment is repeated $10000$ times so that the sizes and powers of the tests can be estimated very accurately.  The type I error $\alpha=0.05$ is fixed in our simulation study.

When $\rho=0$, the null hypothesis in \eqref{approach} is true. The estimated sizes for these test statistics are reported in Table~\ref{table1}.
When $\rho=0.02$, the alternative hypothesis in \eqref{approach} is true, and this indicates a weak dependence among the coordinates of a normal random vector. The estimated powers for these test statistics are given in Table~\ref{table2}.

In terms of the estimated size, a test is considered to be preferable if its estimated size is close to the nominal level ($\alpha=0.05$
in our study).  Table~\ref{table1} indicates that $t_{np}^*$ and $T_{np}^*$ are comparable in terms of the estimated size for tests and the
normal approximation yields significantly larger sizes than the nominal level for both $t_{np}^*$ and $T_{np}^*$ when the dimension $p$ is relatively small.   The test statistics $t_{np}^c$ and $T_{np}^c$ have much better performance than their competitors $t_{np}^*$ and $T_{np}^*$  when $p$ is small as the chi-square approximation is used to determine the corresponding rejection regions. When $p$ is large, the four test
statistics are comparable.

The estimated powers for the four test statistics  are recorded in Table~\ref{table2}. From the table, both $t_{np}^*$ and $T_{np}^*$ have slightly larger powers than $t_{np}^c$ and $T_{np}^c$ for small $p$. This is not surprising since the normal approximation to $t_{np}$ and $T_{np}$
sacrifices the accuracy in the size of the tests when $p$ is small. The performances  of the four test statistics are
similar when $p$ is large.

In summary, we can conclude that the test statistics $t_{np}^c$ and $T_{np}^c$ are consistently accurate in terms of the size
 over the whole range of $p$ and achieve satisfactory power compared with Schott's $t_{np}$ and Mao's $T_{np}$ .
Our simulation study suggests that the normal approximation to Schott's $t_{np}$ and Mao's $T_{np}$ are inferior to the chi-square approximation to  $t_{np}^c$ and $T_{np}^c$ when
$p$ is small. When $p$ is large, the four test statistics under consideration are quite similar in terms of powers and accuracy
in the size. Our simulation also confirms the theoretical consistency in using the normal approximation to
 both $t_{np}$ and $T_{np}$ under the complete independence when $p_n\to\infty$ as $n\to\infty$ regardless of how fast $p_n$ increases with $n$.

%%%%%%%%%%%%%%%%%%%%%%%%%%%%%%%%%%

\begin{table}[t]
\small
%\footnotesize
%\null\vspace{-60pt}
\caption{Size of tests ($\rho=0$)}
\label{table1}
\centering
\begin{tabular}{lrrrrrrrr}
 \hline
  Test Statistic& $n\backslash p$&~~~~~~3&~~~~~ 10&~~~~~ 20& ~~~~~50& ~~~~100 &~~~~200\\
  \hline
  $t_{np}^*$ &$15$ & 0.0718& 0.0583& 0.0614& 0.0563 & 0.0574 & 0.0576\\
           &  $30$ & 0.0725& 0.0611& 0.0598& 0.0510 & 0.0560 & 0.0550\\
           &  $60$ & 0.0711& 0.0593& 0.0559& 0.0544 & 0.0519 & 0.0580\\
           & $100$ & 0.0738& 0.0611& 0.0587& 0.0539 & 0.0542 & 0.0493\\
           & $200$ & 0.0712& 0.0606& 0.0551& 0.0554 & 0.0506 & 0.0497\\
    \cline{2-8}
  $T_{np}^*$ &  $15$ &0.0633 &0.0619 & 0.0605& 0.0554 & 0.0573 & 0.0578\\
           &  $30$ &0.0696 &0.0600 & 0.0599& 0.0522 & 0.0553 & 0.0541\\
           &  $60$ &0.0700 &0.0617 & 0.0568& 0.0556 & 0.0519 & 0.0575\\
           & $100$ &0.0726 &0.0623 & 0.0597& 0.0537 & 0.0551 & 0.0490\\
           & $200$ &0.0713 &0.0606 & 0.0561& 0.0555 & 0.0505 & 0.0494\\
    \cline{2-8}
$t_{np}^c$    &  $15$ &0.0478 &0.0479 & 0.0556& 0.0537 & 0.0566 & 0.0574\\
           &  $30$ &0.0497 &0.0509 & 0.0535& 0.0490 & 0.0550 & 0.0546\\
           &  $60$ &0.0501 &0.0512 & 0.0505& 0.0516 & 0.0507 & 0.0577\\
           & $100$ &0.0526 &0.0512 & 0.0545& 0.0525 & 0.0530 & 0.0490\\
           & $200$ &0.0514 &0.0512 & 0.0494& 0.0531 & 0.0493 & 0.0486\\
   \cline{2-8}
$T_{np}^c$    &  $15$ & 0.0466 &0.0539 & 0.0562& 0.0536 & 0.0558 & 0.0577\\
           &  $30$ & 0.0485 &0.0509 & 0.0544& 0.0508 & 0.0546 & 0.0537\\
           &  $60$ & 0.0503 &0.0515 & 0.0505& 0.0537 & 0.0510 & 0.0572\\
           & $100$ & 0.0520 &0.0524 & 0.0550& 0.0521 & 0.0535 & 0.0487\\
           & $200$ & 0.0519 &0.0512 & 0.0495& 0.0531 & 0.0493 & 0.0490\\
   \hline
\end{tabular}
\end{table}

\begin{table}[t]
%\footnotesize
\small
%\null\vspace{-60pt}
\caption{Power of tests:  $\rho=0.02$}
\label{table2}\centering
\begin{tabular}{lrrrrrrrr}
 \hline
  Test Statistic& $n\backslash p$&~~~~~~3&~~~~~ 10&~~~~~ 20& ~~~~~50& ~~~~100 &~~~~200\\
  \hline
  $t_{np}^*$ &  $15$ & 0.0725 & 0.0656& 0.0647& 0.0765& 0.1003& 0.1557\\
           &  $30$ & 0.0717 & 0.0693& 0.0757& 0.1002& 0.1598& 0.3130\\
           &  $60$ & 0.0811 & 0.0805& 0.0932& 0.1667& 0.3206& 0.6505\\
           & $100$ & 0.0812 & 0.0902& 0.1297& 0.2651& 0.5714& 0.9176\\
           & $200$ & 0.0901 &0.1255 & 0.2175& 0.5834& 0.9413& 0.9996\\
    \cline{2-8}
  $T_{np}^*$ &  $15$ & 0.0641 & 0.0673& 0.0661& 0.0744& 0.1017& 0.1505\\
           &  $30$ & 0.0689 & 0.0715& 0.0756& 0.0984& 0.1583& 0.3096\\
           &  $60$ & 0.0788 & 0.0820& 0.0923& 0.1667& 0.3199& 0.6511\\
           & $100$ & 0.0793 & 0.0909& 0.1303& 0.2646& 0.5706& 0.9173\\
           & $200$ & 0.0894 & 0.1258& 0.2183& 0.5838& 0.9415& 0.9996\\
    \cline{2-8}
$t_{np}^c$    &  $15$ & 0.0494 & 0.0556& 0.0579& 0.0742& 0.0987& 0.1539\\
           &  $30$ & 0.0513 & 0.0581& 0.0685& 0.0977& 0.1574& 0.3116\\
           &  $60$ & 0.0594 & 0.0677& 0.0843& 0.1610& 0.3182& 0.6492\\
           & $100$ & 0.0578 & 0.0777& 0.1212& 0.2594& 0.5674& 0.9166\\
           & $200$ & 0.0649 & 0.1085& 0.2045& 0.5763& 0.9398& 0.9996\\
   \cline{2-8}
$T_{np}^c$    &  $15$ & 0.0466 & 0.0577& 0.0602& 0.0720& 0.1005& 0.1486\\
           &  $30$ & 0.0506 & 0.0604& 0.0686& 0.0954& 0.1566& 0.3083\\
           &  $60$ & 0.0594 & 0.0703& 0.0853& 0.1621& 0.3172& 0.6499\\
           & $100$ & 0.0578 & 0.0778& 0.1211& 0.2589& 0.5672& 0.9164\\
           & $200$ & 0.0647 & 0.1076& 0.2041& 0.5754& 0.9400& 0.9996\\
   \hline
\end{tabular}
\end{table}

%  & $10$ & & & & & & & & \\

\normalsize

\newpage

\section{Proofs}\label{proof}

\noindent{\bf Proof of Theorem~\ref{thm1}.}

We will employ a martingale central limit theorem in
McLeish (1974). Since some details are somewhat similar to those in Mao (2014), we outline our proof as follows.

\textit{Step 1.} Express $r_{ij}$ as $r_{ij}=w_i'w_j$, where $w_1, w_2, \cdots, w_{p_n}$ are independent random vectors that are
uniformly distributed on the surface of the $(n-1)$-sphere. Let $\mathcal{F}_{n\ell}=\sigma(w_1, w_2, \cdots, w_{\ell})$ denote the $\sigma$-algebra generated by $\{w_1, w_2, \cdots, w_{\ell}\}$, see Mao (2014).

\textit{Step 2.} For $2\le \ell\le p_n$\, set $y_{n\ell}=\sigma_{np_n}^{-1}\sum^{\ell-1}_{j=1}\hat{r}_{\ell j}$,
where $\hat{r}_{\ell j}=\frac{r_{\ell j}^2}{1-r_{\ell j}^2}-\frac{1}{n-4}$.
Then, $\{y_{n\ell}, ~\mathcal{F}_{n\ell}, ~2\le \ell\le p_n, ~n\ge 6\}$ form an array of martingale differences, see Mao (2014). Note that $T_{np}^*=\sum^{p_n}_{\ell=2}y_{n\ell}$. According to Theorem 2.3 in McLeish (1974), to show \eqref{s2}, it suffices to prove the following
three conditions:

\noindent (a)  $\displaystyle\sup_{n\ge n_0}E(\max_{2\le \ell\le p_n}(y_{n\ell})^2)<\infty$ for some $n_0$;

\noindent (b)  $\displaystyle\max_{2\le\ell\le p_n}|y_{n\ell}|$ converges to zero in probability;

\noindent (c)  $\displaystyle\sum_{\ell=2}^{p_n}y_{n\ell}^2$ converges to one in probability.

To verify the above three conditions, we need to show that
\begin{equation}\label{2conditions}
\sum^{p_n}_{\ell=2}E(y_{n\ell}^4)\to 0 ~~\mbox{ and } ~~ E(\sum_{\ell=2}^{p_n}y_{n\ell}^2-1)^2\to 0
\end{equation}
as $n\to\infty$.  The second limit implies condition (c) immediately. The first limit implies condition (a),
since
\[
E(\max_{2\le \ell\le p_n}(y_{n\ell})^2)\le \sqrt{E(\max_{2\le \ell\le p_n}y_{n\ell}^4)}\le \sqrt{\sum^{p_n}_{\ell=2}E(y_{n\ell}^4)}\to 0.
\]
Condition (b) follows from the above equation by using the Markov inequality.

It has been proved in Mao (2014) that
\begin{equation}\label{moment}
E(\prod^4_{i=1}\hat{r}_{\ell j_i})=
\left\{
  \begin{array}{ll}
    \frac{12(n-3)(5n^2-27n+40)}{(n-4)^2(n-6)(n-8)(n-10)}, & \hbox{if $j_1=j_2=j_3=j_4$;}\\
    \frac{4(n-3)^2}{(n-4)^2(n-6)^2}, & \hbox{if $\{j_1,j_2,j_3,j_4\}$ forms two distinct pairs;} \\
    0, & \hbox{otherwise.}
  \end{array}
\right.
\end{equation}
Note that $\sigma_{np_n}\sim \frac{p_n}{n}$ as $n\to\infty$. Then, we have
\[
E(y_{n\ell}^4)=\sigma_{np_n}^{-4}\sum_{1\le j_1, j_2, j_3, j_4\le \ell-1}E(\hat{r}_{\ell j_1}\hat{r}_{\ell j_2}\hat{r}_{\ell j_3}\hat{r}_{\ell j_4})
=\sigma_{np_n}^{-4}O\left(\frac{\ell^2}{n^4}\right)
=O\left(\frac{\ell^2}{p_n^4}\right)
\]
uniformly over $2\le \ell\le p_n$ as $n\to\infty$. Therefore, $\sum^{p_n}_{\ell=2}E(y_{n\ell}^4)=O(1/p)\to 0$ as $n\to\infty$. This proves
the first limit in \eqref{2conditions}. Mao (2014) has shown that $E(\sum^{p}_{\ell=2}y_{n\ell}^2)=1$ and
$\sum_{2\le i\ne j\le p_n}E(y_{ni}^2y_{nj}^2)-1=-\frac{2\sigma_{np}^{-4}(n-3)^2p(p-1)(2p-1)}{3(n-4)^4(n-6)^2}$ which is of order $p_n^{-1}$.
Therefore, we have as $n\to\infty$ that
\[
E(\sum_{\ell=2}^{p_n}y_{n\ell}^2-1)^2=E(\sum_{\ell=2}^{p_n}y_{n\ell}^4)+\sum_{2\le i\ne j\le p_n}E(y_{ni}^2y_{nj}^2)-1=O(p_n^{-1})\to 0,
\]
 which yields the second limit in \eqref{2conditions}.   This completes the proof of the theorem. \hfill$\blacksquare$\\

\noindent{\bf Proof of Theorem~\ref{unified}.}
To prove \eqref{chisq}, it suffices to show that for every sequence of integers $\{n_i, ~~i\ge 1\}$, there exists its subsequence
$\{n_{i(j)}, ~j\ge 1\}$ such that \eqref{chisq} holds along $\{n_{i(j)}\}$. Here we choose the subsequence so that $p_{n_{i(j)}}$ converges as
$j\to\infty$. Since $p_{n_{i(j)}}$'s are integers, the limit of $n_{i(j)}$ is a finite integer $p$ or infinity. Therefore, we need to show that
\eqref{chisq} holds along any subsequence of integers $n_i$ such that $p_{n_i}$ is a fixed integer $p$ for all large $i$ or $p_{n_i}\to\infty$
as $i\to\infty$. Since the proof of \eqref{chisq} along a subsequence is the same as the that along the entire sequence, for simplicity,
we will show \eqref{chisq} under the following conditions:
\begin{equation}\label{fixed}
p_n=p\ge 2\mbox{ is a fixed integer for all large }n;
\end{equation}
\begin{equation}\label{infinity}
p_n\to\infty~~~\mbox {as }n\to\infty.
\end{equation}

First, we will show under \eqref{approach} and \eqref{fixed} that
\[
(n-4)T_{np}\td \chi^2_{p(p-1)/2}~~~\mbox{ as }n\to\infty,
\]
which implies \eqref{chisq} since $T_{np}^c=(1+o(1))(n-4)T_{np}+o(1)$.

Express $w_j=z_j/(z_j'z_j)^{1/2}$ for $1\le j\le p$, where $z_j=(z_{j1}, \cdots, z_{j(n-1)})'$, $1\le j\le p$ are i.i.d. random vectors with
$N_{n-1}(0, \mathbf{I}_{n-1})$ distribution.   Write $s_{i,j}=z_i'z_j=\sum^{n-1}_{k=1}z_{ik}z_{jk}$.
By using the multivariate central limit theorem,
\begin{equation}\label{chisq-all}
\frac{1}{\sqrt{n-1}}(s_{2,1}, s_{3,1}, s_{3,2}, \cdots, s_{p,1}, \cdots, s_{p,(p-1)})'\td N_{p(p-1)/2}(0, \mathbf{I}_{p(p-1)/2})
\end{equation}
as $n\to\infty$, which implies that
$\frac{1}{n-1}(s_{2,1}^2, s_{3,1}^2, s_{3,2}^2, \cdots, s_{p,1}^2, \cdots, s_{p,(p-1)}^2)'$ converges in distribution to a random vector
whose $p(p-1)/2$ components are independent random variables having a chi-squared distribution with $1$ degree of freedom.
By the law of large numbers,  $\frac{z_i'z_i'}{n-1}=1+o_p(1)$ for $i=1, \cdots, p$, which implies
\[
\max_{1\le i\le p}\left|\frac{z_i'z_i}{n-1}-1\right|=o_p(1)~~~\mbox{ as }n\to\infty.
\]
Therefore,
 \begin{equation}\label{rij}
r_{ij}^2=s_{i,j}^2/((z_i'z_i)(z_j'z_j))=\frac{s_{i,j}^2}{(n-1)^2}(1+o_p(1)),
\end{equation}
 it follows that
\[
\frac{(n-4)r_{ij}^2}{1-r_{ij}^2}=\frac{s_{i,j}^2}{n-1}(1+o_p(1)),
\]
which implies that
\[
(n-4)T_{np}=\frac{\displaystyle\sum_{1\le j<i\le p}s_{ij}^2}{n-1}(1+o_p(1))\td \chi^2_{p(p-1)/2} ~~\mbox{ as }n\to\infty.
\]

Now assume \eqref{infinity}  and the null hypothesis in \eqref{approach} hold. In this case, we can apply Theorem~\ref{thm1} directly. It follows from
\eqref{TC} and \eqref{s2} that
\[
\frac{T_{np}^c-\frac{p(p-1)}{2}}{\sqrt{p(p-1)}}=T_{np}^*
\overset{d}\to N(0,1),
\]
 which implies that
\begin{equation}\label{ineq1}
\sup_x\left|P\left(\frac{T_{np}^c-\frac{p(p-1)}{2}}{\sqrt{p(p-1)}}\le x\right)-\Phi(x)\right|\to 0~~~\mbox{ as } n\to\infty,
\end{equation}
where $\Phi(x)$ is the standard normal cumulative distribution function. Also, notice that a chi-squared random variable
with $p(p-1)/2$ degrees of freedom can be written as the sum of $p(p-1)/2$ independent and identically distributed random variables
having a chi-squared distribution with $1$ degree of freedom. From the classic central limit theorem, we have
\[
\frac{\chi^2_{p(p-1)/2}-\frac{p(p-1)}{2}}{\sqrt{p(p-1)}}\td N(0,1),
\]
 and thus
\begin{equation}\label{ineq2}
\sup_x\left|P\left(\frac{\chi^2_{p(p-1)/2}-\frac{p(p-1)}{2}}{\sqrt{p(p-1)}}\le x\right)-\Phi(x)\right|\to 0.
\end{equation}
Therefore, by combining \eqref{ineq1} and \eqref{ineq2} and using the triangle inequality we have
\begin{eqnarray*}
&&\sup_x|P(T_{np}^c\le x)-P(\chi^2_{p(p-1)/2}\le x)|\\
&=&\sup_x\left|P\left(\frac{T_{np}^c-\frac{p(p-1)}{2}}{\sqrt{p(p-1)}}\le x\right) -P\left(\frac{\chi^2_{p(p-1)/2}-\frac{p(p-1)}{2}}{\sqrt{p(p-1)}}\le x\right)\right|\\
&\le &\sup_x\left|P\left(\frac{T_{np}^c-\frac{p(p-1)}{2}}{\sqrt{p(p-1)}}\le x\right)-\Phi(x)\right|+\sup_x \left|P\left(\frac{\chi^2_{p(p-1)/2}-\frac{p(p-1)}{2}}{\sqrt{p(p-1)}}\le x\right)-\Phi(x)\right|\\
&\to& 0
\end{eqnarray*}
as $n\to\infty$. This completes the proof of \eqref{chisq}. \hfill$\blacksquare$

\noindent{\bf Proof of Theorem~\ref{schott-all}.}

We will sketch the proof. We continue to use the notation in the proof of Theorem~\ref{unified}. As in the proof of Theorem~\ref{thm1}, write
\begin{equation}\label{freeway}
r_{ij}=w_i'w_j, ~~1\le i, j\le p.
\end{equation}

\noindent (i) First, we need to show \eqref{schott}, i.e.,
\begin{equation}\label{Tn}
t_{np_n}^*=\frac{\sum^{p_n-1}_{i=2}\sum^{i-1}_{j=1}r_{ij}^2-\frac{p_n(p_n-1)}{2(n-1)}}{\tau_{np_n}}\td N(0,1),
\end{equation}
under the assumption that $p_n\to\infty$ as $n\to\infty$.

Set
\[
z_{n\ell}=\sum^{\ell-1}_{i=1}r_{\ell i}^2-\frac{\ell-1}{n-1}.
\]
Then, $\{z_{n\ell}, ~\mathcal{F}_{n\ell}, ~ 2\le \ell\le p_n,~n\ge 1\}$ form an array of martingale differences. See, e.g., Schott (2005).
It suffices to show that
\begin{equation}\label{CLT}
\frac{\sum^{p_n}_{\ell=2}z_{n\ell}}{\tau_{np_n}}\td N(0,1).
\end{equation}

We will use martingale approach like that in Schott (2005). In view of Corollary 3.1 in Hall and Heyde (1980), the martingale central limit theorem \eqref{CLT} holds if the following two conditions hold:
\begin{equation}\label{martingale1}
\frac{1}{\tau_{np_n}^2}\sum^{p_n}_{\ell=2}E(z_{n\ell}^2I(|z_{n\ell}|\ge \varepsilon \tau_{np_n})|\mathcal{F}_{n(\ell-1)})\to 0~~~\mbox{ in probability }
\end{equation}
for every $\varepsilon>0$, and
\begin{equation}\label{martingale2}
\frac{1}{\tau_{np_n}^2}\sum^{p_n}_{\ell=2}E(z_{n\ell}^2|\mathcal{F}_{n(\ell-1)})\to 1~~~\mbox{ in probability}.
\end{equation}

It has been shown in Schott (2005), pp. 955 that
\[
E\left(\sum^{p_n}_{\ell=2}E(z_{n\ell}^2|\mathcal{F}_{n(\ell-1)})\right)=\tau_{np_n}^2.
\]
Thus, we have
\begin{eqnarray*}
\Delta_n:&=&\sum^{p_n}_{\ell=2}E(y_{n\ell}^2|\mathcal{F}_{n, \ell-1})-\tau_{np}^2\\
&=&\frac{2}{(n-1)(n+1)}\sum^{p_n}_{\ell=2}\sum^{\ell-1}_{i=1}\sum^{\ell-1}_{j=1,j\ne i}\left(r_{ij}^2-\frac{1}{n-1}\right)\\
&=&\frac{2}{(n-1)(n+1)}\sum^{p_n}_{\ell=3}\sum_{1\le i\ne j\le \ell-1}\left(r_{ij}^2-\frac{1}{n-1}\right)\\
&=&\frac{4}{(n-1)(n+1)}\sum^{p_n}_{\ell=3}\sum_{1\le j<i\le \ell-1}\left(r_{ij}^2-\frac{1}{n-1}\right)\\
&=&\frac{4}{(n-1)(n+1)}\sum_{1\le j<i\le p_n-1}(p_n-i)\left(r_{ij}^2-\frac{1}{n-1}\right).
\end{eqnarray*}
It is easy to verify from Schott (2005) that
\[
E\left(r_{ij}^2-\frac{1}{n-1}\right)\left(r_{st}^2-\frac{1}{n-1}\right)
=\left\{
 \begin{array}{ll}
  0, & \hbox{ if } (i,j)\ne (s,t),\\
 \frac{3}{(n-1)(n+1)}-\frac{1}{(n-1)^2}=\frac{2n-4}{(n-1)(n+1)}, & \hbox{ if } (i,j)\ne (s,t).
  \end{array}
   \right.
\]
Then, we have
\begin{eqnarray*}
E(\Delta_n^2)&=&\frac{16}{(n-1)^2(n+1)}\sum_{1\le j<i\le p_n-1}(p_n-i)^2\frac{2(n-2)}{(n-1)^2(n+1)}\\
&=&\frac{32(n-2)}{(n-1)^4(n+1)^3}\sum_{1<i\le p_n-1}(p_n-i)^2(i-1)\\
&=&O\left(\frac{p_n^4}{n^6}\right),
\end{eqnarray*}
which implies that
\begin{equation}\label{lyapunov}
\frac{E(\Delta_n^2)}{\tau_{np_n}^4}=O\left(\frac{1}{n^2}\right)\to 0 ~\mbox{ as }n\to\infty.
\end{equation}

Next, we verify that
\begin{equation}\label{4thmoment}
\frac{1}{\tau_{np}^4}\sum^{p_n}_{\ell=2}E(z_{n\ell}^4)=o(1)~~~\mbox{as }n\to\infty.
\end{equation}

Set $q_{\ell i}=r_{\ell i}^2-\frac{1}{n-1}$. Then
\[
z_{n\ell}=\sum^{\ell-1}_{i=1}q_{\ell i}.
\]
Note that $r_{\ell i}=w_{\ell}'w_{i}$. Conditional on $w_{\ell}$, $r_{\ell 1}, \cdots, r_{\ell (\ell-1)}$ are i.i.d.
Set $c_r=E(r_{\ell 1}^{2r}|w_{\ell})$, $1\le r\le 4$. Then
\[
c_1=\frac{1}{n-1}, ~c_2=\frac{3}{(n-1)(n+1)},~c_3=\frac{15}{(n-1)(n+1)(n+3)},
\]
and
\[
c_4=\frac{105}{(n-1)(n+1)(n+3)(n+5)}.
\]

Set $d_r=E(q_{\ell 1}^r|w_{\ell})=E((r_{\ell 1}^2-\frac{1}{n-1})^r|w_{\ell})$. Then
\[
d_1=0, ~~~~d_2=c_2-\left(\frac{1}{n-1}\right)^2=\frac{2(n-2)}{(n-1)^2(n+1)},
\]
\[
d_3=c_3-3c_2\frac{1}{n-1}+3c_1\frac{1}{(n-1)^2}-\frac{1}{(n-1)^3}=O\left(\frac{1}{n^3}\right),
\]
and
\[
d_4=c_4-4c_3\frac{1}{n-1}+6c_2\frac{1}{(n-1)^2}-4c_1\frac{1}{(n-1)^4}+\frac{1}{(n-1)^4}=O\left(\frac{1}{n^4}\right).
\]
Since
\[
E(z_{n\ell}^4)=E(E(z_{n\ell}^4|w_{\ell}))=E(E((\sum^{\ell-1}_{i=1}q_{\ell i})^4|w_{\ell}))=(\ell-1)d_4+6(\ell-1)(\ell-2)d_2^2
\]
for $2\le \ell\le p$, we obtain
\[
\sum^{p_n}_{\ell=2}E(z_{n\ell}^4)=O\left(\frac{p_n^3}{n^4}\right).
\]
Then, it follows that
\[
\frac{1}{\tau_{np}^4}\sum^{p_n}_{\ell=2}E(z_{n\ell}^4)=O\left(\frac{1}{p_n}\right)\to 0~~\mbox{ as }n\to\infty,
\]
which implies \eqref{4thmoment}.  \eqref{martingale1} and \eqref{martingale2} can be easily verified from
\eqref{lyapunov} and \eqref{4thmoment}.   Therefore, we obtain \eqref{CLT}.

\noindent{(ii)} For the proof of \eqref{t-chisq}  we can use the arguments in the proof of Theorem~\ref{unified}.
First, under assumption \eqref{fixed}, we have from \eqref{chisq-all} and \eqref{rij} that $(n-1)t_{np}\td \chi^2_{p(p-1)/2}$
as $n\to\infty$.  The rest of the proof follows exactly the same lines as that in the proof of Theorem~\ref{unified}
by using \eqref{schott}. The details are omitted. \hfill$\blacksquare$

\vspace{20pt}

\noindent{\bf Acknowledgements:}  
%We would like to thank two reviewers for their constructive suggestions that have led to improvement in the layout and readability of the %paper.
Chang's research was supported in part by the Major Research Plan of the National Natural Science Foundation of China (91430108), the National Basic Research Program (2012CB955804), the National Natural Science Foundation of China (11171251), and the Major Program of Tianjin University of Finance and Economics (ZD1302).

\end{document}